\documentclass{amsart}
\usepackage{amsmath,amsthm}
\usepackage{graphics}
\usepackage{color}
\usepackage{epsfig}

\graphicspath{{./Figures/}}

\theoremstyle{plain}

\newtheorem*{thm_main}{Theorem}
\newtheorem*{corblank}{Corollary}

\newcommand{\comment}[1]{}

\newcommand{\Z}{\ensuremath{\mathbb{Z}}}

\newcommand{\bdry}{\ensuremath{\partial}}

\begin{document}

\title[Surgery descriptions and volumes of Berge knots II]{Surgery descriptions and volumes of Berge knots II: Descriptions on the minimally twisted five chain link}

\author{Kenneth L. Baker}
\address{Department of Mathematics, University of Georgia \\ Athens, Georgia 30602}
\email{kb@math.uga.edu}

\thanks{This work was partially supported by a graduate traineeship from the VIGRE Award at the University of Texas at Austin and a VIGRE postdoc under NSF grant number DMS-0089927 to the University of Georgia at Athens.}



\begin{abstract}
Using Kirby Calculus, we explicitly pass from Berge's R-R descriptions of ten families of knots with lens space surgeries to surgery descriptions on the minimally twisted five chain link (MT5C).  Since the MT5C admits a strong involution, we also give the corresponding tangle descriptions.  
\end{abstract}

\maketitle

\section{Introduction}

In \cite{berge:skwsyls} Berge describes twelve families of knots that admit lens space surgeries.  These knots are referred to as {\em Berge knots} and appear to comprise all knots in $S^3$ known to have lens space surgeries.  Ten of these families, (I)-(VI) being the knots in solid tori with surgeries yielding solid tori and (IX)-(XII) being the `sporadic' knots, are described via R-R diagrams.  

\begin{thm_main}\label{thm:main}
The knots in Berge's families (I)-(VI) and (IX)-(XII) admit surgery descriptions on the minimally twisted five chain link.
\end{thm_main}
\begin{figure}[h]
\centering
\input{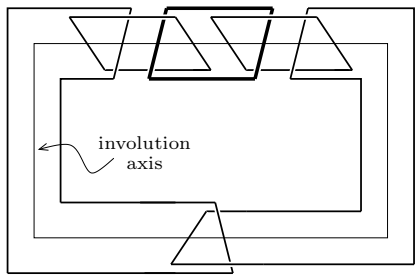}
\caption{The minimally twisted five chain link, MT5C, with an axis of strong involution.}
\label{fig:MT5C}
\end{figure}
Figure~\ref{fig:MT5C} shows the minimally twisted five chain link, MT5C for short.  The proof of the theorem is given in \S\ref{sec:proofofmain}.  We use Kirby Calculus (see e.g.\ \cite{rolfsen} or \cite{gompfstipsicz:4makc}) to pass from these R-R diagrams to surgery descriptions on the MT5C.  Since the MT5C is strongly invertible, we give the corresponding tangle descriptions of these surgeries in \S\ref{sec:tangles}.

As an immediate consequence of Thurston's Hyperbolic Dehn Surgery Theorem~\cite{thurston:gt3m} and the fact that the MT5C is a hyperbolic link (e.g.\ Theorem 5.1 (ii) of \cite{neumannreid}), our theorem has the following corollary.

\begin{corblank}
Volumes of hyperbolic knots in Berge's families (I)-(VI) and (IX)-(XII) are bounded above by the volume of the MT5C.
\end{corblank} 

Families (VII) and (VIII) are the knots which lie as essential simple closed curves on the fiber of the trefoil and figure eight knot respectively.  In the prequel~\cite{baker:sdavobkI} we show that they contain hyperbolic knots of arbitrarily large volume.  Consequentially, they cannot all be all be described by surgery on the MT5C.  Nevertheless, each of them admits a surgery description on some minimally twisted chain link. 

\subsection{Acknowledgements}
The author wishes to thank both John Luecke for his direction and many useful conversations and Yuichi Yamada for his comments.

\section{Proof of Theorem}\label{sec:proofofmain}

\begin{figure}[h!]
\centering
\input{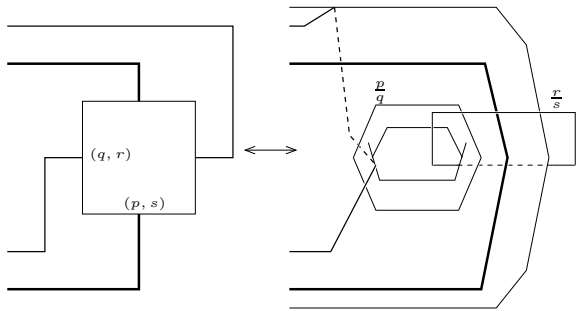}
\caption{Passing from a R-R diagrams to a surgery description.}
\label{fig:Diagram-to-sfce}
\end{figure}

We prove the theorem by exhibiting a passage from Berge's R-R diagrams to surgery descriptions on the MT5C, or its reflection.  Please refer to \cite{berge:skwsyls} for Berge's original descriptions of these knots and the notation conventions.  The first step in this passage translates the R-R diagram to a knot on a genus $2$ Heegaard surface for $S^3$ together with surgery instructions on a surrounding link.  Half of the correspondence is shown in Figure~\ref{fig:Diagram-to-sfce}.  The other half is obtained by rotating the pictures $180^\circ$ in the plane of the page.

The first six families of Berge knots arise from knots in solid tori with surgeries yielding solid tori.  Where relevant we take $n, p, q, r, s, K \in \Z$, $|ps-qr| = 1$ and $\epsilon = \pm 1$.

\subsection{Type (I), torus knots}

\begin{figure}[h!]
\centering
\input{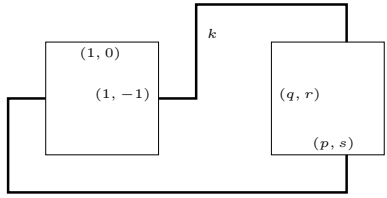}
\caption{Berge's R-R diagram for type (I) knots (torus knots).} 
\label{fig:I-Diagram}
\end{figure}

\begin{figure}[h!]
\centering
\input{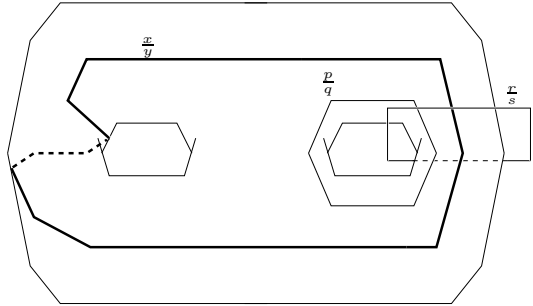}
\caption{Type (I) knots on Heegaard surface via surgeries.} 
\label{fig:I-sfce}
\end{figure}

We pass from Berge's diagram (Figure~\ref{fig:I-Diagram})  to a corresponding realization of these knots $k$ on a genus 2 Heegaard surface via surgeries (Figure~\ref{fig:I-sfce}).  In Figure~\ref{fig:I-sfce}, $\frac{x}{y} = \frac{1}{0}$ corresponds to the meridional filling of $k$ while $\frac{x}{y} = -\frac{1}{1}$ corresponds to the lens space surgery of $k$.  

\begin{figure}[h!]
\centering
\input{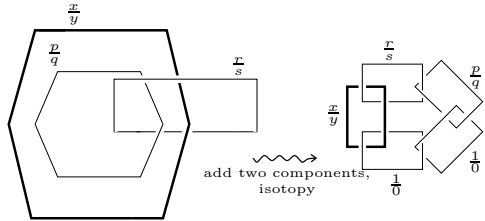}
\caption{Type (I) knots.  Isotopies and Kirby Calculus moves ending in the MT5C.} 
\label{fig:I-isotopies}
\end{figure}

After dropping the Heegaard surface from the picture,we add components with meridional framings, perform isotopies, and do Kirby Calculus as shown in Figure~\ref{fig:I-isotopies} to get a description of $k$ as surgery on the MT5C.

\subsection{Type (II), cables about torus knots}

\begin{figure}[h!]
\centering
\input{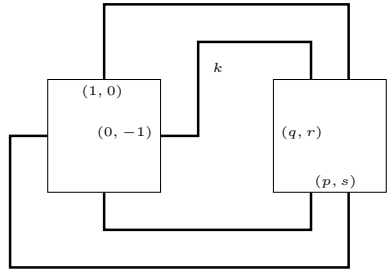}
\caption{Berge's R-R diagram for type (II) knots (cables about torus knots).} 
\label{fig:II-Diagram}
\end{figure}

\begin{figure}[h!]
\centering
\input{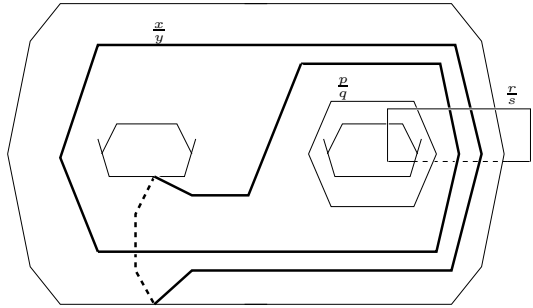}
\caption{Type (II) knots on Heegaard surface via surgeries.} 
\label{fig:II-sfce}
\end{figure}

We pass from Berge's diagram (Figure~\ref{fig:II-Diagram}) to a corresponding realization of these knots $k$ on a genus 2 Heegaard surface via surgeries (Figure~\ref{fig:II-sfce}).  In Figure~\ref{fig:II-sfce}, $\frac{x}{y} = \frac{1}{0}$ corresponds to the meridional filling of $k$ while $\frac{x}{y} = -\frac{1}{1}$ corresponds to the lens space surgery of $k$.  

\begin{figure}[h!]
\centering
\input{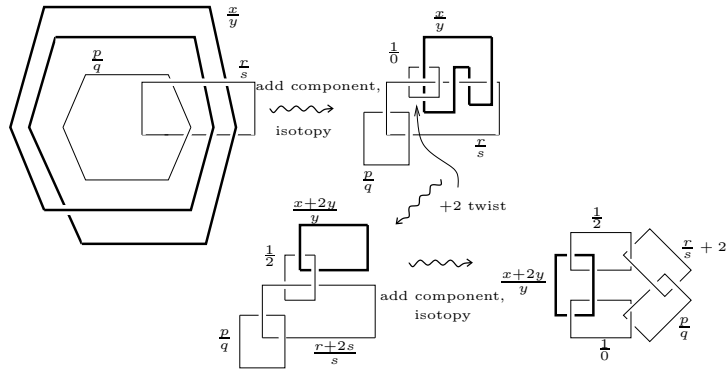}
\caption{Type (II) knots.  Isotopies and Kirby Calculus moves ending in the MT5C.} 
\label{fig:II-isotopies}
\end{figure}

After dropping the Heegaard surface from the picture, we add components with meridional framings, perform isotopies, and do Kirby Calculus as shown in Figure~\ref{fig:II-isotopies} to get a description of $k$ as surgery on the MT5C.  Notice that the surgery of $+\frac{1}{1}$ on the component corresponding to $k$ is the lens space surgery.  The trivial (i.e. $S^3$) surgery is still $\frac{1}{0}$.

\subsection{Type (III)}

\begin{figure}[h!]
\centering
\input{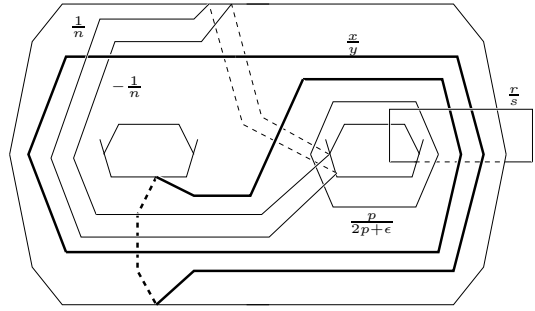}
\caption{Type (III) knots on Heegaard surface via surgeries.} 
\label{fig:III-sfce}
\end{figure}

We pass from Berge's diagram (Figure~\ref{fig:III-Diagram}) to a corresponding realization of these knots $k$ on a genus 2 Heegaard surface via surgeries (Figure~\ref{fig:III-sfce}).  In Figure~\ref{fig:III-sfce}, $\frac{x}{y} = \frac{1}{0}$ corresponds to the meridional filling of $k$ while $\frac{x}{y} = -\frac{1}{1}$ corresponds to the lens space surgery of $k$.  Also, $\frac{r}{s} = \frac{2 \epsilon + (2p + \epsilon) K}{\epsilon + p K}$. 

\begin{figure}[h!]
\centering
\input{Figures/III-Diagram.pstex_t}
\caption{Berge's R-R diagram for type (III) knots.} 
\label{fig:III-Diagram}
\end{figure}

\begin{figure}[h!]
\centering
\input{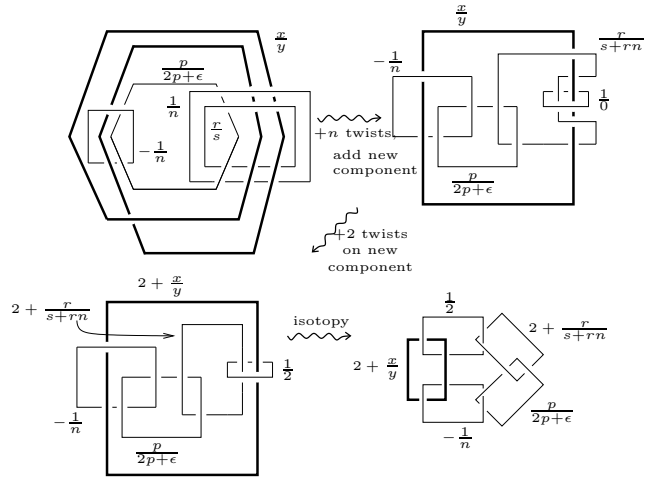}
\caption{Type (III) knots.  Isotopies and Kirby Calculus moves ending in the MT5C.} 
\label{fig:III-reduced}
\end{figure}

We drop the Heegaard surface to get the first link of Figure~\ref{fig:III-reduced}.  The surgeries on the two parallel components in the first link of Figure~\ref{fig:III-reduced} can be amalgamated into a surgery on a single component.  From here we perform isotopies and Kirby Calculus to arrive at the MT5C.   Notice that the surgery of $+\frac{1}{1}$ on the component corresponding to $k$ on the MT5C is the lens space surgery.  The trivial (i.e. $S^3$) surgery is still $\frac{1}{0}$.

\subsection{Type (IV)}

We pass from Berge's diagram (Figure~\ref{fig:IV-Diagram})  to a corresponding realization of these knots $k$ on a genus 2 Heegaard surface via surgeries (Figure~\ref{fig:IV-sfce}).  In Figure~\ref{fig:IV-sfce}, $\frac{x}{y} = \frac{1}{0}$ corresponds to the meridional filling of $k$ while $\frac{x}{y} = -\frac{1}{1}$ corresponds to the lens space surgery of $k$.  Also, $\frac{r}{s} = \frac{2 \epsilon + (2p + \epsilon) K}{\epsilon +p K}$.

\begin{figure}[h!]
\centering
\input{Figures/IV-Diagram.pstex_t}
\caption{Berge's R-R diagram for type (IV) knots.} 
\label{fig:IV-Diagram}
\end{figure}

\begin{figure}[h!]
\centering
\input{Figures/IV-sfce.pstex_t}
\caption{Type (IV) knots on Heegaard surface via surgeries.} 
\label{fig:IV-sfce}
\end{figure}

After dropping the Heegaard surface from the picture, we add components with meridional framings, perform isotopies, and do Kirby Calculus as shown in Figure~\ref{fig:IV-isotopies} to get a description of $k$ as surgery on the MT5C.  Notice that the surgery of $-\frac{1}{1}$ on the component corresponding to $k$ is the lens space surgery.  The trivial (i.e. $S^3$) surgery is still $\frac{1}{0}$.

\begin{figure}[h!]
\centering
\input{Figures/IV-isotopies-newer.pstex_t}
\caption{Type (IV) knots.  Isotopies and Kirby Calculus moves ending in the MT5C.} 
\label{fig:IV-isotopies}
\end{figure}



\subsection{Type (V)}

\begin{figure}[h!]
\centering
\input{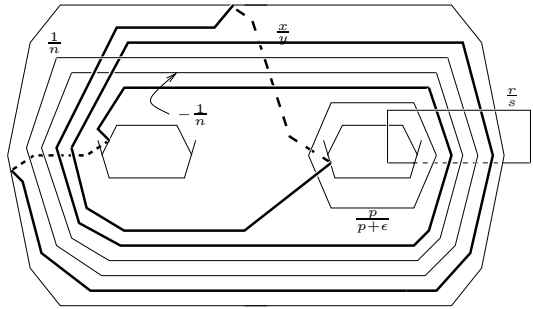}
\caption{Type (V) knots on Heegaard surface via surgeries.} 
\label{fig:V-sfce}
\end{figure}

We pass from Berge's diagram (Figure~\ref{fig:V-Diagram})  to a corresponding realization of these knots $k$ on a genus 2 Heegaard surface via surgeries (Figure~\ref{fig:V-sfce}).  In Figure~\ref{fig:V-sfce}, $\frac{x}{y} = \frac{1}{0}$ corresponds to the meridional filling of $k$ while $\frac{x}{y} = -\frac{5}{1}$ corresponds to the lens space surgery of $k$.  Also, $\frac{r}{s} = \frac{2 \epsilon + (2p + \epsilon) K}{\epsilon +p K}$.

\begin{figure}[h!]
\centering
\input{Figures/V-Diagram.pstex_t}
\caption{Berge's R-R diagram for type (V) knots.} 
\label{fig:V-Diagram}
\end{figure}

\begin{figure}[h!]
\centering
\input{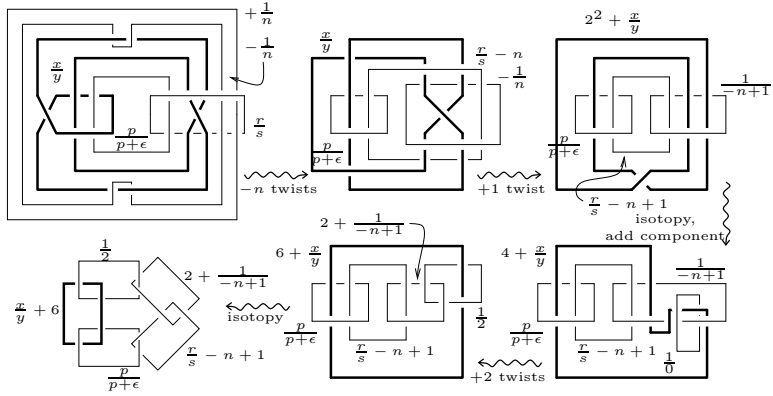}
\caption{Type (V) knots.  Isotopies and Kirby Calculus moves ending in the MT5C.} 
\label{fig:V-isotopies}
\end{figure}

After dropping the Heegaard surface from the picture, we amalgamate two parallel components and do Kirby Calculus as shown in Figure~\ref{fig:V-isotopies} to get a description of $k$ as surgery on the MT5C.  Notice that the surgery of $+\frac{1}{1}$ on the component corresponding to $k$ is the lens space surgery.  The trivial (i.e. $S^3$) surgery is still $\frac{1}{0}$.

\subsection{Type (VI)}

\begin{figure}[h!]
\centering
\input{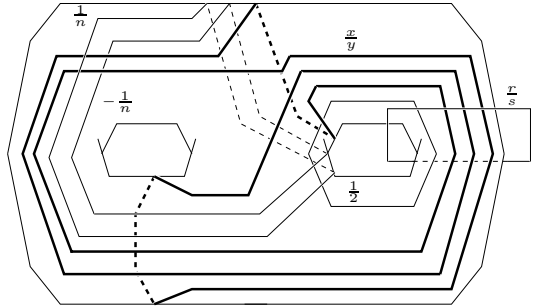}
\caption{Type (VI) knots on Heegaard surface via surgeries.} 
\label{fig:VI-sfce2}
\end{figure}

We pass from Berge's diagram (Figure~\ref{fig:VI-Diagram})  to a corresponding realization of these knots $k$ on a genus 2 Heegaard surface via surgeries (Figure~\ref{fig:VI-sfce2}).  In Figure~\ref{fig:VI-sfce2}, $\frac{x}{y} = \frac{1}{0}$ corresponds to the meridional filling of $k$ while $\frac{x}{y} = -\frac{5}{1}$ corresponds to the lens space surgery of $k$. 

\begin{figure}[h!]
\centering
\input{Figures/VI-Diagram2.pstex_t}
\caption{Berge's R-R diagram for type (VI) knots.} 
\label{fig:VI-Diagram}
\end{figure}

\begin{figure}[h!]
\centering
\input{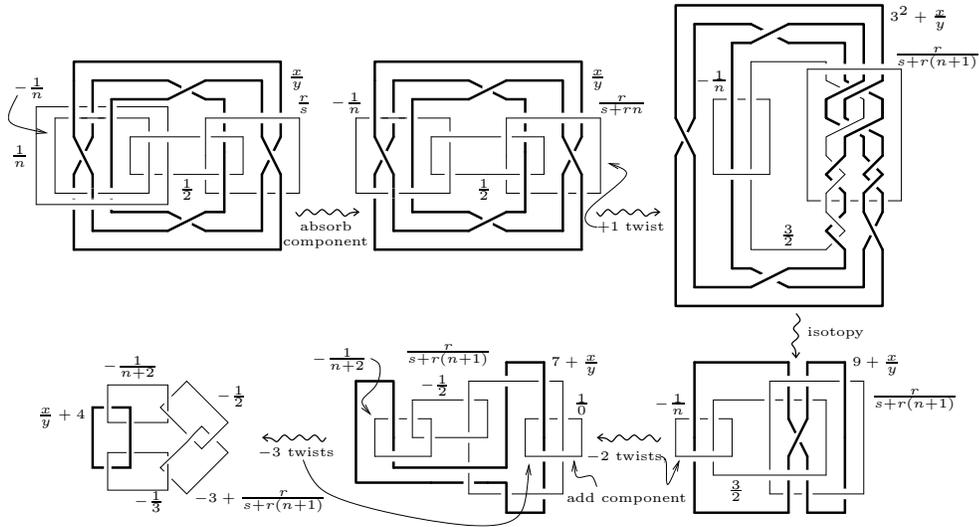}
\caption{Type (VI) knots.  Isotopies and Kirby Calculus moves ending in the MT5C.} 
\label{fig:VI-isotopies}
\end{figure}

After dropping the Heegaard surface from the picture, we amalgamate the surgeries on two parallel components in the first step of Figure~\ref{fig:VI-isotopies}.  We then continue with Kirby Calculus to get a description of $k$ as surgery on the MT5C.  Notice that the surgery of $-\frac{1}{1}$ on the component corresponding to $k$ is the lens space surgery.  The trivial (i.e. $S^3$) surgery is still $\frac{1}{0}$.

\subsection{Sporadic knots type a) and b)}

\begin{figure}[h!]
\centering
\input{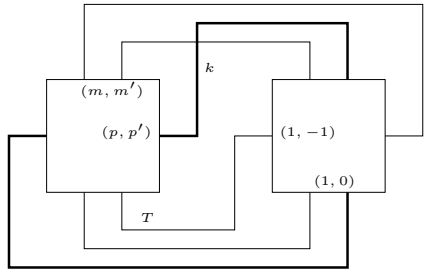}
\caption{Berge's sporadic knot types a) and b) diagram.}
\label{fig:SpAB-Diagram}
\end{figure}

\begin{figure}[h!]
\centering
\input{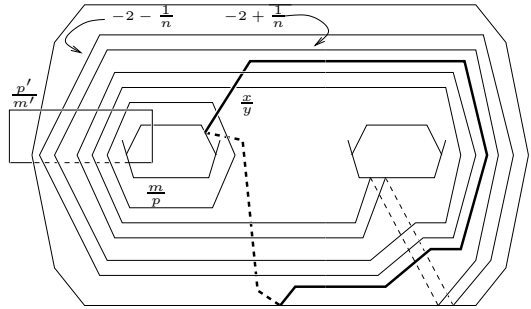}
\caption{Sporadic knot types a) and b) on Heegaard surface via surgeries.}
\label{fig:SpAB-sfce}
\end{figure}

We pass from Berge's diagram (Figure~\ref{fig:SpAB-Diagram})  to a corresponding realization of these knots $k$ on a genus 2 Heegaard surface via surgeries (Figure~\ref{fig:SpAB-sfce}).  In Figure~\ref{fig:SpAB-sfce}, $\frac{x}{y} = \frac{1}{0}$ corresponds to the meridional filling of $k$ while $\frac{x}{y} = \frac{0}{1}$ corresponds to the lens space surgery of $k$.  Here $n \in Z$.  For type a) knots $(p,p',m,m') = (1,1,2,3)$.  For type b) knots $(p,p',m,m') = (2,1,3,2)$.

\begin{figure}[h!]
\centering
\input{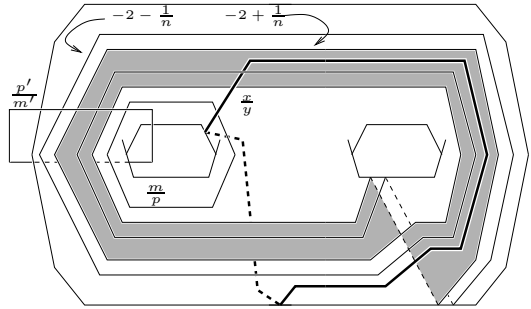}
\caption{Sporadic types a) and b). One component bounds a M\"obius band.}
\label{fig:SpAB-sfce-annulus}
\end{figure}

As shown in Figure~\ref{fig:SpAB-sfce-annulus}, one component bounds a M\"obius band. The framing on that component induced by the M\"obius band agrees with the framing from the surface.

\begin{figure}[h!]
\centering
\input{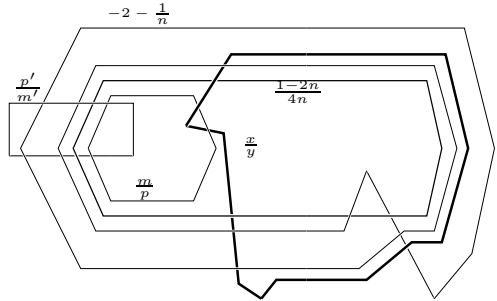}
\caption{Sporadic types a) and b).  Core curve of M\"obius band replaces the boundary of M\"obius band.}
\label{fig:SpAB-reduced}
\end{figure}

\begin{figure}[h!]
\centering
\input{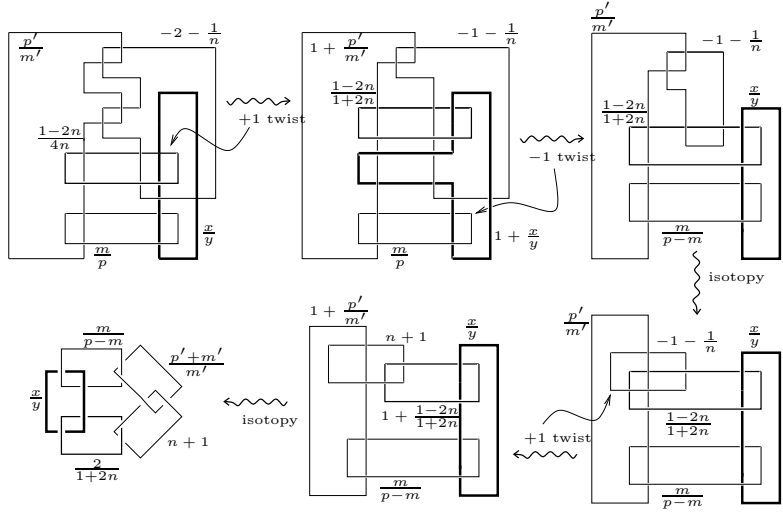}
\caption{Sporadic types a) and b).  Isotopies and Kirby Calculus moves ending in the MT5C.}
\label{fig:SpAB-isotopies}
\end{figure}

In Figure~\ref{fig:SpAB-reduced}, we replace the component that bounds the M\"obius band with the core curve of the M\"obius band and the corresponding surgery instructions.  We continue, after an isotopy, with Kirby Calculus in Figure~\ref{fig:SpAB-isotopies} to get a description of $k$ as surgery on the MT5C. 

\subsection{Sporadic knots type c) and d)}

\begin{figure}[h!]
\centering
\input{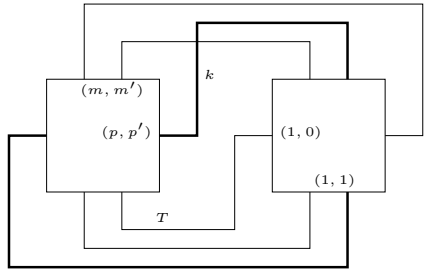}
\caption{Berge's sporadic knot types c) and d) diagram.}
\label{fig:SpCD-Diagram}
\end{figure}

We pass from Berge's diagram (Figure~\ref{fig:SpCD-Diagram})  to a corresponding realization of these knots $k$ on a genus 2 Heegaard surface via surgeries (Figure~\ref{fig:SpCD-sfce}).  In Figure~\ref{fig:SpCD-sfce}, $\frac{x}{y} = \frac{1}{0}$ corresponds to the meridional filling of $k$ while $\frac{x}{y} = +\frac{1}{1}$ corresponds to the lens space surgery of $k$.  Here $n \in Z$.  For type c) knots $(p,p',m,m') = (4,-3,3,-2)$.  For type d) knots $(p,p',m,m') = (3,-5,2,-3)$.

\begin{figure}[h!]
\centering
\input{Figures/SpCD-sfce.pstex_t}
\caption{Sporadic knot types c) and d) on Heegaard surface via surgeries.}
\label{fig:SpCD-sfce}
\end{figure}

\begin{figure}[h!]
\centering
\input{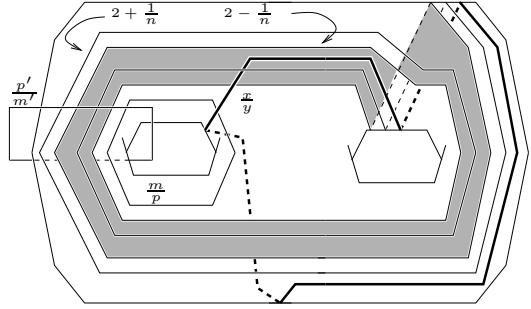}
\caption{Sporadic types c) and d). One component bounds a M\"obius band.}
\label{fig:SpCD-sfce-annulus}
\end{figure}

As shown in Figure~\ref{fig:SpCD-sfce-annulus},  one component bounds a M\"obius band. The framing on that component induced by the M\"obius band agrees with the framing from the surface.

\begin{figure}[h!]
\centering
\input{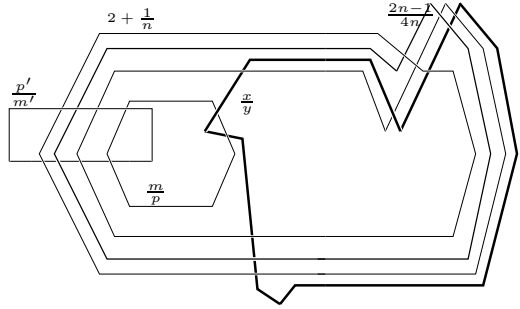}
\caption{Sporadic types c) and d).  Core curve of M\"obius band replaces the boundary of M\"obius band.}
\label{fig:SpCD-reduced}
\end{figure}

\begin{figure}[h!]
\centering
\input{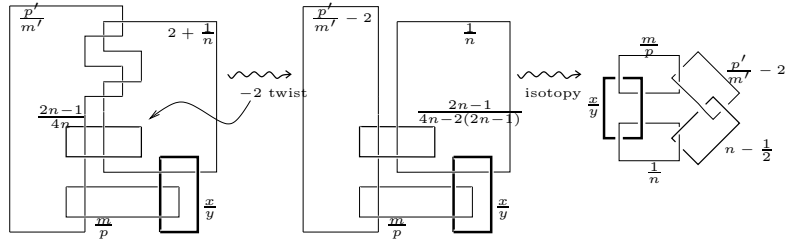}
\caption{Sporadic types c) and d).  Isotopies and Kirby Calculus moves ending in the MT5C.}
\label{fig:SpCD-isotopies}
\end{figure}

In Figure~\ref{fig:SpCD-reduced}, we replace the component that bounds the M\"obius band with the core curve of the M\"obius band and the corresponding surgery instructions.  We continue, after an isotopy, with Kirby Calculus in Figure~\ref{fig:SpCD-isotopies} to get a description of $k$ as surgery on the MT5C.

\subsection{Summary of Berge knots as surgeries on the MT5C.}\label{subsec:summary}

\begin{figure}[h]
\centering
\input{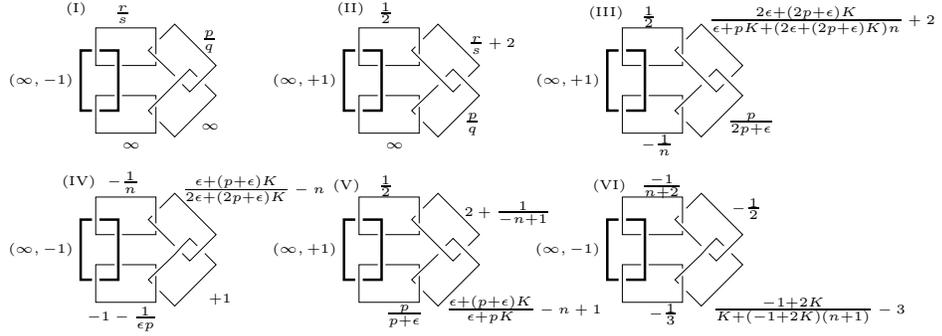}
\caption{Berge knot types (I)---(VI)}
\label{fig:MT5C-BergeI-VI}
\end{figure}

\begin{figure}[h]
\centering
\input{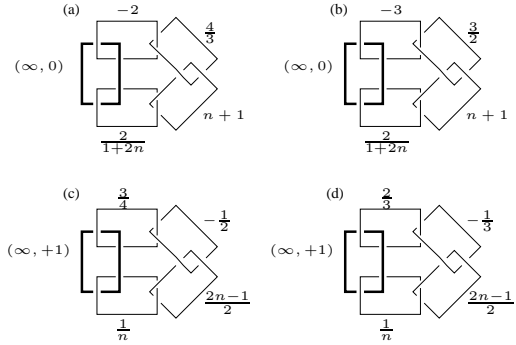}
\caption{Sporadic Berge knot types a)---d)}
\label{fig:MT5C-BergeSp}
\end{figure}

We conclude the proof of the theorem by summarizing the results above.
Figure~\ref{fig:MT5C-BergeI-VI} shows Berge knot types (I) --- (VI) as surgeries on the minimally twisted five chain link or its mirror.  Figure~\ref{fig:MT5C-BergeSp} shows Berge sporadic knot types a) --- d) as surgeries on the the minimally twisted five chain link or its mirror.  The component corresponding to the Berge knot is shown with a pair of surgeries 
$(\rho_{S^3}, \rho_{\textrm{LensSp}})$ 
where $\rho_{S^3}$ is the surgery slope yielding $S^3$ and $\rho_{\textrm{LensSp}}$ is the surgery slope yielding the lens space.  Recall also that $n, p, q, r, s, \in \Z$, $|ps-qr|=1$ and $\epsilon = \pm1$.

\qed

\section{Tangles}\label{sec:tangles}
Because the MT5C is strongly invertible, we may rephrase the minimally twisted chain link surgery descriptions of Berge's families of knots as tangle descriptions.  We include these tangle descriptions as they may be of interest to others and provide an alternate means of verification of the above surgery descriptions.  

On the remaining boundary component of each of the the tangles below, a choice of rational tangles to be inserted is indicated.  Inserting the first rational tangle (which is the $\infty$ tangle in each case) renders the tangle into the unknot.  As the double branched cover of the unknot in $S^3$ is $S^3$ again, this corresponds to the trivial surgery on the Berge knot.  Inserting the second rational tangle yields a two-bridge link, a link that may be decomposed into two rational tangles.  As the double branched cover of a two-bridge knot is a lens space, this corresponds to the lens space surgery on the Berge knot.

\subsection{Conventions}\label{subsec:tangleconv}
 A {\em tangle} $(B,t)$ is a pair consisting of a punctured $3$-sphere $B$ and a properly embedded collection $t$ of disjoint arcs and simple closed curves.  Two tangles $(B_1, t_1)$ and $(B_2, t_2)$ are homeomorphic if there is a homeomorphism of pairs
\[h \colon (B_1, t_1) \to (B_2, t_2). \]

A boundary component $(\bdry B_0, t \cap \bdry B_0)$ of a tangle $(B,t)$ is a sphere $\bdry B_0$ together with some finite set of points $p=t \cap \bdry B_0$.  Here we consider the situation where $p$ consists of just four points.

Given a sphere $S$ and set $p$ of four distinct points on $S$,
 a {\em framing} of $(S, p)$ is an ordered pair of (unoriented) simple closed curves $(\hat{m}, \hat{l})$ on $S - N(p)$ such that each curve separates different pairs of points of $p$.
 
\begin{figure}
\centering
\input{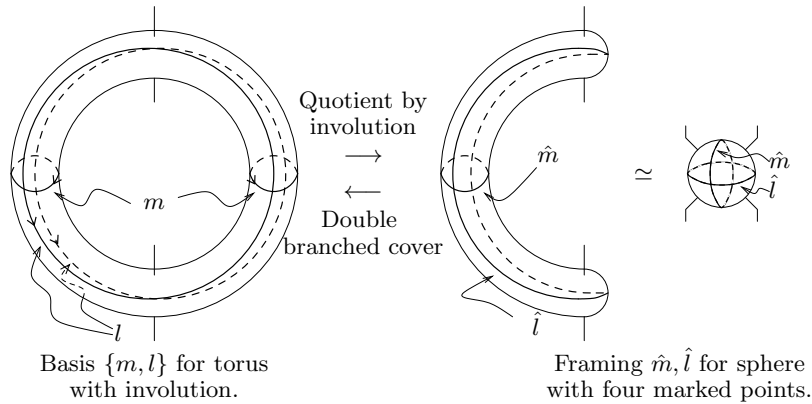}
\caption{The correspondence of bases and framings.}
\label{fig:tangleframing}
\end{figure}

 Let $(S, p)$ be a sphere with four points with framing $(\hat{m}, \hat{l})$.  The double cover of $S$ branched over $p$ is a torus.  Single components, say $m$ and $l$, of the lifts of the framing curves $\hat{m}$ and $\hat{l}$ when oriented so that $m \cdot l = +1$ (with respect to the orientation of the torus) form a basis for the torus.  Similarly, a basis on a torus induces a framing on the sphere with four points obtained by the quotient of an involution that fixes four points on the torus.  See Figure~\ref{fig:tangleframing}.
  
This implies the correspondence between inserting into a boundary component of a tangle the rational tangle $p/q$ and the $p/q$ Dehn filling of a torus boundary component in the double branched cover of the tangle.

\begin{figure}[h]
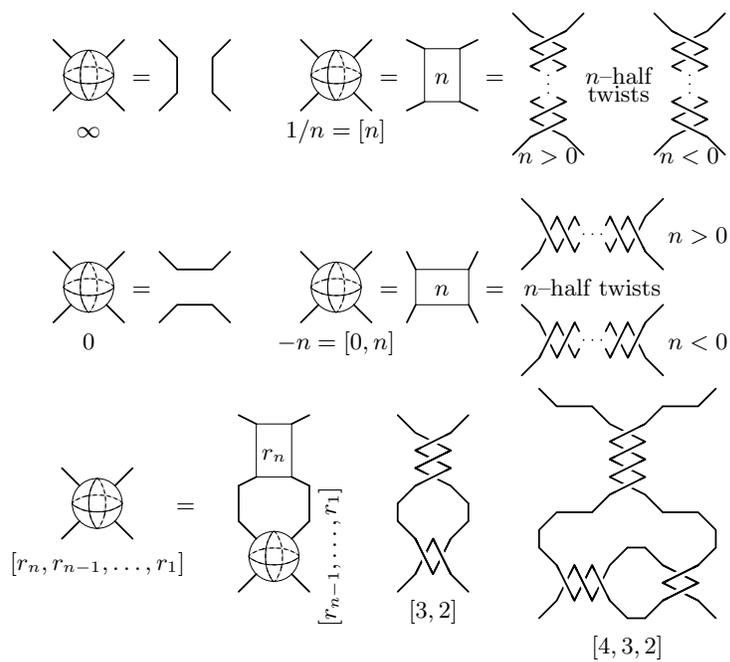

\[\begin{array}{c}
\input{Figures/tanglelegend1.pstex_t}\\
\input{Figures/tanglelegend2.pstex_t}
\end{array}\]
\caption{Tangle legend}
\label{fig:tanglelegend}
\end{figure}

Figure~\ref{fig:tanglelegend} shows by example the conventions we use for rational tangles.  Here, $[a,b, \dots, c]$ denotes the continued fraction $1/(a - 1/(b- 1/( \dots -1/c)))$.  Notice that $1 = [1] = [0,-1]$. 


\comment{
\begin{figure}[h]
\centering
\input{Figures/tanglelegend1.pstex_t}
\caption{Tangle legend}
\label{fig:tanglelegend1}
\end{figure}

\begin{figure}[h]
\centering
\input{Figures/tanglelegend2.pstex_t}
\caption{Tangle legend and example}
\label{fig:tanglelegend2}
\end{figure}
}

\subsection{Tangle descriptions}
Please refer to \cite{baker:sdavobkI} for an illustration of the passage between a surgery description on the MT5C to a tangle description.  We now present the tangle descriptions obtained from the surgery descriptions in \S\ref{subsec:summary} above using the tangle conventions in \S\ref{subsec:tangleconv}.  Note that some of the surgery descriptions in \S\ref{subsec:summary} are on the mirror of the MT5C and these tangle descriptions reflect that accordingly.  In each tangle below,  where presented with the choice $(\infty, \delta)$ for $\delta \in\{-1,0,+1\}$ inserting the rational tangle $\infty$ yields the unknot and inserting the rational tangle $\delta$ yields a $2$--bridge link.

\begin{figure}[h!]
\centering
\input{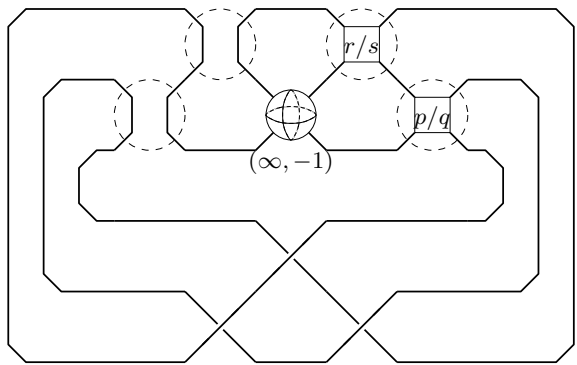}
\caption{Tangle description of Berge knot type (I), torus knots.}
\label{fig:tangleBergeI}
\end{figure}

\begin{figure}[h!]
\centering
\input{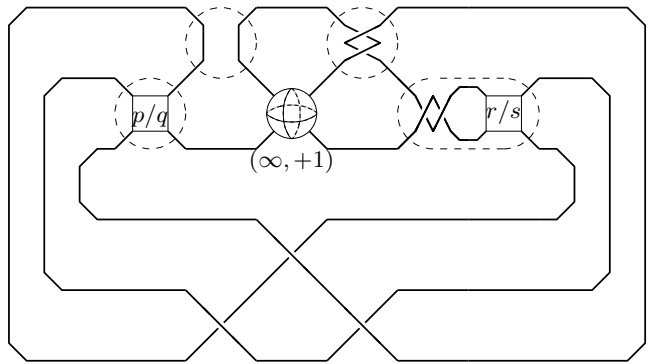}
\caption{Tangle description of Berge knot type (II), cables about torus knots.}
\label{fig:tangleBergeII}
\end{figure}

\begin{figure}[h!]
\centering
\input{Figures/5chaintangleIII.pstex_t}
\caption{Tangle description of Berge knot type (III).}
\label{fig:tangleBergeIII}
\end{figure}

\begin{figure}[h!]
\centering
\input{Figures/5chaintangleIV.pstex_t}
\caption{Tangle description of Berge knot type (IV).}
\label{fig:tangleBergeIV}
\end{figure}

\begin{figure}[h!]
\centering
\input{Figures/5chaintangleV.pstex_t}
\caption{Tangle description of Berge knot type (V).}
\label{fig:tangleBergeV}
\end{figure}

\begin{figure}[h!]
\centering
\input{Figures/5chaintangleVI.pstex_t}
\caption{Tangle description of Berge knot type (VI).}
\label{fig:tangleBergeVI}
\end{figure}

\begin{figure}[h!]
\centering
\input{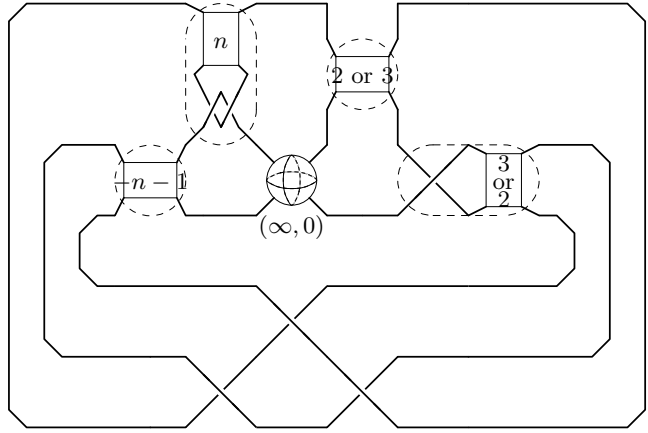}
\caption{Tangle description of Berge knot types (IX) and (X), Sporadics a) and b).}
\label{fig:tangleBergeAB}
\end{figure}

\begin{figure}[h!]
\centering
\input{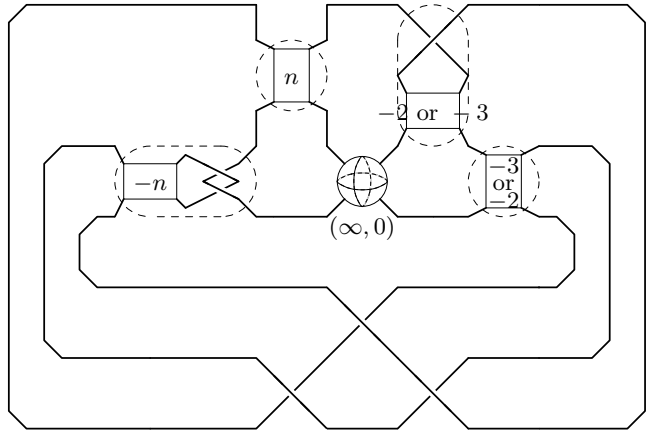}
\caption{Tangle description of Berge knot type (XI) and (XII), Sporadics c) and d).}
\label{fig:tangleBergeCD}
\end{figure}

\clearpage
\bibliography{MathBiblio}
\bibliographystyle{hamsplain}
\end{document}